\documentclass[12pt]{article}
\usepackage{graphicx}
\usepackage{amsmath}
\usepackage[dvips]{epsfig}
\usepackage{amssymb}
\def\thesection{\arabic{section}.}
\makeatletter
\renewcommand{\section}{\@startsection
  {section}%
  {2}%
  {0mm}%
  {\baselineskip}%
  {0.5 \baselineskip}%
  {\centering}}
\makeatother
\begin{document}

\title { A new approach to  $q$-zeta function }
\author{ Taekyun Kim   \\\\
Institute of Science Education,\\
 Kongju National University,
Kongju 314-701,  Korea\\  e-mail: tkim@kongju.ac.kr
  }

\date{}
\maketitle

 {\footnotesize {\bf Abstract}\hspace{1mm}
We construct the new $q$-extension of Bernoulli numbers and
polynomials in this paper. Finally we consider the $q$-zeta
functions which interpolate the new $q$-extension of Bernoulli
numbers and polynomials.



\bigskip
\section{Introduction }
\bigskip

Throughout this paper  $\mathbb{Z},  \mathbb{Z}_p,  \mathbb{Q}_p$
and  $\mathbb{C}_p$  will be denoted by  the ring of rational
integers, the ring of $p$-adic  integers, the field of $p$-adic
rational numbers and the
  completion of algebraic closure of   $\mathbb{Q}_p$,
  respectively,
cf. [7, 8, 9, 10].  Let $\nu_p$ be the normalized exponential
valuation of $\mathbb{C}_p$ with $|p|_p=p^{-\nu_p(p)}=p^{-1}.$
When one talks of $q$-extension, $q$ is variously considered as an
indeterminate, a complex number $q\in \mathbb{C},$ or $p$-adic
number $q \in \Bbb C_p .$  If $q\in  \mathbb{C}_p,$ then we
normally assume  $|q-1|_p<p^{-\frac{1}{p-1}},$ so that $q^x=\exp(x
\log q)$ for $|x|_p\leq 1.$
 If $q\in \Bbb C$, then we normally assume that
$|q|<1.$ For $f\in UD( \mathbb{Z}_p,\mathbb{C}_p )=\{ f | f :
\mathbb{Z}_p \rightarrow \mathbb{C}_p \mbox{ is uniformly
differentiable function} \}, $ the $p$-adic $q$-integral (or
$q$-Volkenborn integration) was defined as
$$I_q(f)=\int_{ \mathbb{Z}_p}f(x)d\mu_q(x)=\lim_{N\rightarrow \infty}
\frac{1}{[p^N]_q}\sum_{x=0}^{p^N-1}f(x)q^x, \eqno(1) $$ where
$[x]_q =\dfrac{1-q^x}{1-q},$   cf. [1, 2, 3, 4, 11]. Thus we note
that
$$I_1(f)=\lim_{q \rightarrow 1} I_q(f)=\int_{ \mathbb{Z}_p}f(x)d\mu_1(x)=\lim_{N\rightarrow \infty}
\frac{1}{p^N } \sum_{0\leq x<p^N}f(x), \mbox{ cf. [4, 11]}.
\eqno(2)$$ By (2), we easily see that
$$I_1(f_1)=I_1(f)+f^\prime (0), \mbox{ cf. [5, 6, 7]},  \eqno(3)$$
where $f_1(x)=f(x+1),  f^\prime (0)=
\dfrac{d}{dx}f(x)\big|_{x=0}.$

  In [8], the $q$-Bernoulli polynomials are defined by
$$\beta_{n}^{(h)}(x,q) =\int_{ \mathbb{Z}_p}[x+x_1]_q^n q^{x_1(h-1)}
d\mu_q(x_1), \mbox{ for } h \in \mathbb{Z}. \eqno(4)$$ In this
paper we consider the new $q$-extension of Bernoulli numbers and
polynomials. The main purpose of this paper is to construct the
new $q$-extension of zeta function and $L$-function which
interpolate the above new $q$-extension of  Bernoulli numbers at
negative integer.

\bigskip
\section{ On the New $q$-Extensions of Bernoulli numbers and polynomials }
\bigskip
In (3), if we take $f(x)=q^{hx}e^{xt}, $ then we have
$$\int_{ \mathbb{Z}_p}q^{hx}e^{xt} d\mu_1(x)=\dfrac{h \log q +t}{q^h e^t -1}. \eqno(5)$$
for $| t | \leq p^{-\frac{1}{p-1}}, h \in \mathbb{Z}$.

Let us define the $(h,q)$-extension of Bernoulli numbers and
polynomials as follows:
$$F_q^{(h)}(t)=\dfrac{h \log q +t}{q^h e^t -1}= \sum_{n=0}^{\infty} B_{n,q}^{(h)} \dfrac{t^n}{n!},$$
$$F_q^{(h)}(t,x)=\dfrac{h \log q +t}{q^h e^t -1}e^{xt}=
 \sum_{n=0}^{\infty} B_{n,q}^{(h)}(x) \dfrac{t^n}{n!}.\eqno(6)$$
Note that $B_{n,q}^{(h)}(0) = B_{n,q}^{(h)},$  $  \lim_{q
\rightarrow 1}B_{n,q}^{(h)}=B_n, $ where $B_n$ are the $n$-th
Bernoulli numbers. By (5) and (6), we obtain the following Witt's
formula.

\bigskip
{ \bf Theorem 1.} For $h \in \mathbb{Z},$ $ q \in \mathbb{C}_p $
with $| 1-q |_p  \leq p^{-\frac{1}{p-1}}$, we have
$$\int_{ \mathbb{Z}_p}q^{hx}x^n d\mu_1(x)= B_{n,q}^{(h)}, $$
$$ \int_{ \mathbb{Z}_p}q^{hy}(x+y)^n d\mu_1(y)= B_{n,q}^{(h)}(x). \eqno(7)$$
\bigskip
By above theorem, we easily see that
$$ B_{n,q}^{(h)}(x)= \sum_{k=0}^n \binom nk x^{n-k} B_{k,q}^{(h)}. \eqno(8)$$

 Let $d$ be any fixed positive integer  with $(p,d)=1$.
Then we set
 $$\aligned
&X=X_d=\varprojlim_N ( \mathbb{Z}/dp^N \mathbb{Z}), X_1=\mathbb{Z}_p,\\
&X^*=\bigcup_{ 0<a<dp} a+dp \mathbb{Z}_p,\\
&a+dp^N \mathbb{Z}_p=\{x\in X \mid x\equiv a\pmod{dp^N}\},
\endaligned$$
where $a\in \mathbb{ Z}$ with $0\leq a<dp^N$. Note that $$\int_{
\mathbb{Z}_p} f(x) d\mu_1(x)=\int_{ X} f(x) d\mu_1(x) ,$$ for
$f\in UD( \mathbb{Z}_p,\mathbb{C}_p )$, cf. [1, 2, 3, 4, 10].

In Eq. (7), it is easy to see that
$$\aligned
& B_{k,q}^{(h)}(x) =\int_X (x+t)^k q^{ht}d\mu_1(t) = \lim_{l
\rightarrow
\infty} \dfrac{1}{mp^l}\sum_{n=0}^{mp^l-1}q^{hn}(x+n)^k \\
&=\frac{1}{m} \lim_{l \rightarrow \infty}
\dfrac{1}{p^l}\sum_{i=0}^{m-1}
\sum_{n=0}^{p^l-1}q^{h(i+mn)}(x+i+mn)^k  =m^{k-1}
\sum_{i=0}^{m-1}q^{hi} B_{k,q^m}^{(h)}\left(
\dfrac{x+i}{m}\right).
\endaligned $$
Therefore we have the below theorem.

\bigskip
{ \bf Theorem 2.}  For any positive integer $m$, we have
$$B_{k,q}^{(h)}(x)
= m^{k-1} \sum_{i=0}^{m-1}q^{hi} B_{k,q^m}^{(h)} \left(
\dfrac{x+i}{m} \right), \mbox{ for } k \geq 0.$$
\bigskip

Let $\chi$ be the Dirichlet  character with conductor $d \in
\mathbb{Z_+}$. Then we define the $(h,q)$-extension of generalized
Bernoulli numbers attached to $\chi$.  For $n \geq 0$, define
$$B_{n,q, \chi }^{(h)}
= \int_{ X} \chi (x) q^{hx} x^n  d\mu_1(x). \eqno(9) $$
 By (9), we easily see that

$$\aligned
B_{n,q, \chi }^{(h)} &  = \lim_{l \rightarrow \infty}
\dfrac{1}{dp^l}\sum_{x=0}^{dp^l-1}\chi (x) q^{hx}x^n =\frac{1}{d}
\lim_{l \rightarrow \infty} \dfrac{1}{p^l}\sum_{i=0}^{d-1}
\sum_{x=0}^{p^l-1}\chi(i+dx) q^{h(i+dx)}(i+dx)^n \\
&=\frac{1}{d} \sum_{i=0}^{d-1} \chi(i) q^{hi}\int_{\Bbb Z_p }
q^{hdx}(i+dx)^n d\mu_1(x)=d^{n-1} \sum_{i=0}^{d-1}\chi(i)q^{hi}
B_{n,q^d}^{(h)}\left( \dfrac{i}{d}\right).
\endaligned $$
Therefore we obtain the below lemma.

\bigskip
{ \bf Lemma 3.}  For  $d \in \mathbb{Z_+}$,  we have
$$B_{k,q, \chi }^{(h)}
= d^{k-1} \sum_{i=0}^{d-1}\chi(i)q^{hi} B_{k,q^d}^{(h)}\left(
\dfrac{i}{d}\right), \mbox{ for } n \geq 0. \eqno(10)$$
\bigskip

By induction in Eq.(3), we easily see that
$$I_1(f_b)=I_1(f)+ \sum_{i=0}^{b-1} f^\prime (i),
\text{ where $f_b(x)=f(x+b), b \in \mathbb{Z_+}.$}\eqno(11)$$

In Eq.(11), if  we take $f(x)=q^{hx}e^{tx}\chi(x)$, then we have
$$ I_1(e^{tx}q^{hx} \chi(x))=\dfrac{\sum_{i=0}^{d-1}\left( t e^{it} \chi(i) q^{hi}
+e^{ti}(h \log q) q^{hi} \chi(i) \right)}{q^{hd}e^{dt}-1}.
\eqno(12)$$
By (12) and (9), we can give the generation function
of $ B_{n,q,\chi}^{(h)} $ as follows:
$$F_{q, \chi }^{(h)}(t)   = \dfrac{ \sum_{i=0}^{d-1}\left( t e^{it}
\chi(i) q^{hi} +e^{ti}(h \log q) q^{hi} \chi(i)
\right)}{q^{hd}e^{dt}-1} = \sum_{n=0}^{\infty} B_{n,q, \chi}^{(h)}
\dfrac{t^n}{n!}. \eqno(13) $$

\bigskip
\section{ The analogue of zeta function }
\bigskip

In this section  we assume that $q\in \mathbb{C} $ with $ |q|<1$.
Let $\Gamma(s)$ be the gamma function. By (6), we can readily see
that

$$
\aligned &  \dfrac{1}{\Gamma (s)}\int_{0}^\infty t^{s-2}e^{-t}
F_q^{(h)}(-t) dt =\frac{1}{\Gamma (s)} \int_{0}^\infty
t^{s-2}e^{-t}
\left \{ \dfrac{-t}{q^he^{-t}-1}+ \dfrac{h \log q}{q^h e^{-t}-1} \right \} dt \\
&=\frac{1}{\Gamma (s)} \int_{0}^\infty t^{s-1}e^{-t}
 \dfrac{1}{1-q^he^{-t}}dt -
 \dfrac{h \log q}{\Gamma (s)} \int_{0}^\infty \dfrac{t^{s-2} e^{-t}}{1-q^h e^{-t}}  dt \\
&=\sum_{n=0}^\infty q^{nh} \dfrac{1}{\Gamma (s)} \int_{0}^\infty
t^{s-1}e^{-(n+1)t}dt - h
\log q \sum_{n=0}^\infty q^{nh}  \int_{0}^\infty t^{s-1}e^{-(n+1)t}dt \\
&=\sum_{n=1}^\infty \dfrac{q^{(n-1)h} }{n^s}-\dfrac{ \log
q^h}{s-1}\sum_{n=1}^\infty \dfrac{q^{(n-1)h} }{n^{s-1}}.
\endaligned
\eqno(14)
$$
Using (14), we define the new $q$-extensions of zeta functions as
follows:

\bigskip
{ \bf Definition 4.}   For  $s \in \mathbb{C}, x \in \mathbb{
R^+}$,  we define
$$\zeta_q^{(h)}(s)= \sum_{n=1}^\infty
\dfrac{q^{(n-1)h}}{n^s}- \dfrac{ h \log q}{s-1}\sum_{n=1}^\infty
\dfrac{q^{(n-1)h} }{n^{s-1}},\eqno (14-a)$$

$$\zeta_q^{(h)}(s, x)= \sum_{n=0}^\infty
\dfrac{q^{nh}}{(n+x)^s}- \dfrac{ h \log q}{s-1}\sum_{n=0}^\infty
\dfrac{q^{nh} }{(n+x)^{s-1}}.\eqno (14-b)$$

\bigskip

\bigskip
{ \bf Remark.}  By (14-a) and (14-b), we easily see that
$\zeta_q^{(h)}(s)=\zeta_q^{(h)}(s, 1)$.  Also, we note that
$\zeta_q^{(h)}(s)$ is analytic continuation for $R(s)>1$.
\bigskip

Using Mellin transforms in Eq.(6), we obtain
$$ \dfrac{1}{\Gamma (s)}\int_{0}^\infty t^{s-2}
F_q^{(h)}(-t,x) dt= \zeta_q^{(h)}(s,x). \eqno(15) $$ By (15) and
(6), we readily see that
$$\zeta_q^{(h)}(s,x)= \sum_{n=0}^\infty \dfrac{(-1)^n B_{n,q}^{(h)} (x)}{n!}
  \dfrac{1}{\Gamma (s)}\int_{0}^\infty t^{s-2+n} dt.$$
Therefore we obtain the following:

\bigskip
{ \bf Theorem 5.}  For $n \in \mathbb{N}$, we have
$$ \zeta_q^{(h)}(1-n ,x)= -\dfrac{ B_{n,q}^{(h)} (x)}{n} $$
\bigskip

By Mellin transforms  and  (13), we note that
$$ \dfrac{1}{\Gamma (s)}\int_{0}^\infty F_{q,\chi}^{(h)}(-t) t^{s-2} dt
= \sum_{n=1}^\infty \dfrac{q^{nh}\chi(n)}{n^s}- \dfrac{ h \log
q}{s-1}\sum_{n=1}^\infty \dfrac{q^{nh}\chi(n) }{n^{s-1}}.
\eqno(16)
$$
Thus we can define the new $q$-extension of Dirichlet $L$-function
as follows:

\bigskip
{ \bf Definition 6.}   For  $s \in \mathbb{C}$,  we define
$$L_q^{(h)}(s, \chi)= \sum_{n=1}^\infty \dfrac{q^{nh}\chi(n)}{n^s}- \dfrac{ h \log
q}{s-1}\sum_{n=1}^\infty \dfrac{q^{nh}\chi(n) }{n^{s-1}}.$$
\bigskip

By (16) and (13), we have the following

\bigskip
{ \bf Theorem 7.}   For  $ n \in \mathbb{N}$,  we obtain
$$L_q^{(h)}(1-n, \chi)= - \dfrac{B_{n,q,\chi}^{(h)}}{n}.$$

 \end{document}